\documentclass[11pt]{article}
\usepackage{graphicx}
\usepackage{amssymb,amsmath}
\usepackage{multirow}
\usepackage{setspace}

\newcommand{\bu}{\mathbf{u}}
\newcommand{\bv}{\mathbf{v}}
\newcommand{\EE}{\mathrm{E}}
\newcommand{\sR}{\scriptscriptstyle R}

\newtheorem{theorem}{Theorem}

\title{Continuous and Tractable models for the Variation of Evolutionary Rates}
\author{Thomas Lepage\footnote{Dept. of Mathematics and Statistics, McGill University, lMontr\'{e}al}, Stephan Lawi\footnote{Laboratoire de Probabilit\'{e}s et Mod\`{e}les Al\'{e}atoires, Universit\'{e} Pierre et Marie Curie, Paris}, Paul Tupper\footnotemark[1]\,\, and David Bryant\footnotemark[1]}

\begin{document}

\maketitle

\vfill
\noindent{\bf Corresponding author}\\
David Bryant\\
McGill Centre for Bioinformatics\\
3775 University\\
 Montr\'eal, Qu\'ebec H3A 2B4\\
 Canada.\\
  ph: 1 514 398-4578. fax: 1 514 398-3387. \\
  email: bryant@mcb.mcgill.ca\\
  
  \pagebreak


\begin{abstract}
We propose a continuous model for evolutionary rate variation
across sites and over the tree and derive exact transition
probabilities under this model. Changes in rate are modelled using
the CIR process, a diffusion widely used in financial
applications. The model directly extends the standard gamma
distributed rates across site model, with one additional parameter
governing changes in rate down the tree. The parameters of the
model can be estimated directly from two well-known statistics:
the index of dispersion and the gamma shape parameter of the rates
across sites model. The CIR model can be readily incorporated into
probabilistic models for sequence evolution. We provide here an
exact formula for the likelihood of a three taxa tree. Larger
trees can be evaluated using Monte-Carlo methods.
\end{abstract}
 
\vspace{3cm}

\pagebreak

\section{Introduction}

Understanding evolutionary rates and how they vary is one of the
central concerns of molecular evolution. It has been clearly shown
that inadequate models of rate variation, between lineages and
between loci, can dramatically affect the accuracy of phylogenetic
inference~\cite{Y93, LC02, SI02}. The dependency of molecular
dating on evolutionary rate models is even more critical: we will
only obtain precise divergence time estimates from molecular data
once we can model the rate at which sequences
evolve~\cite{AY02,KT01}.

Modelling the evolutionary rate is made difficult by the number
and variety of factors influencing it. The base rate of mutation
can vary because of changes in the accuracy of transcription
machinery \cite{B86}, DNA repair mechanisms \cite{W89}, and
metabolic rate \cite{MP93}.  At the cellular level, selective
pressures can lead to variation of rate between loci and over
time, as evidenced by differential rates of the three codon
position \cite{GY94,MG94}, the slower evolutionary rate of highly
expressed genes \cite{PPH01}, and the effect of tertiary structure
on patterns of sequence conservation \cite{RJKGT03}.

Selection also affects the evolutionary rate at the level of
populations. For the most part, the only mutations that affect
phylogenetics are those that are fixed in the population. Hence
evolutionary rate is a combination of mutation rate and fixation
rate. Fluctuations in population size, generation times, and
environmental pressures affect fixation rates and thereby
influence evolutionary rate \cite{OT90, T87, L69}.

Because of this complexity, the strategies employed for modelling evolutionary rate have tended to be statistical in nature. As with all statistical inference, there is an iterative sequence of model formulation, model assessment, and model improvement. The aim is to construct a model that accurately explains the observed variation  but is as simple, and tractable, as possible.\\

Our goal in this paper is to derive a continuous model for rate evolution that avoids many of the problems of existing approaches. We base our model on the CIR process, a continuous Markov process that is widely used in finance to model interest rates~\cite{CIR85}. As we shall see, the model fits well into existing protocols for phylogenetic inference. The process has a stationary distribution given by a gamma distribution and yet, unlike the rates-across-sites (RAS) model of Uzzell and Corbin \cite{UC71}, the rate is allowed to vary along lineages. The CIR model adds only one parameter to the RAS model, and this parameter can be estimated directly from the index of dispersion or the autocorrelation (see below). Furthermore, we can derive exact transition probabilities when we incorporate CIR based rate variation into the standard models for sequence evolution. \\

The outline of the paper is as follows:
\begin{itemize}
\item In the following section we summarise the key
characteristics of models for rate evolution, and show how
existing models are classified with respect to these
characteristics. \item In Section~\ref{CIR} we present the CIR
model for rate evolution and discuss its basic properties. \item
In Section~\ref{MutModel} we derive transition probabilities for
standard mutation models where the rate is described as a Markov
process. \item In Section~\ref{Sectionmgf} we focus on the case
where the rate is modelled by a CIR process. \item In
Section~\ref{sec:3taxa} we extend this one step further to derive an
expression for likelihood of a three-taxa tree using a mutation
model with rate determined by the CIR process. We note that
three-taxa trees are often used to study differences in
evolutionary rate.
\end{itemize}
We conclude with an outline of future work and work in progress.

In a companion paper (in preparation) we describe the
incorporation of this model into software for Bayesian
phylogenetic inference, and use this to show how our model
captures important information lost in standard RAS approaches.

\section{Properties of models for rate variation} \label{Prop}

In this section we examine several important characteristics that can be used to distinguish, and choose between, different models for rate variation. We discuss how the different existing models fit into this scheme and summarise the differences between them in Table~\ref{ModelSummary}.

The rate of evolution for a given locus at time $t \geq 0$ is denoted by $R_t$.  For each $t>0$, $R_t$ is a non-negative random variable, and different models of rate evolution give different distributions for the rates $R_t$, $t \geq 0$.

Here and throughout the paper we will restrict out attention to
\emph{Markov processes}. That is, for any $t_1 \leq t_2 \leq t_3$,
we assume that $R(t_3)$ conditioned on $R(t_2)$ is independent of
$R(t_1)$. In other words, the future depends on the past only
through the present.\\

\noindent {\em Property I: Continuous or  Discontinuous Sample Paths}\\

\begin{figure}[h]

\begin{center}
\includegraphics[scale = 0.6]{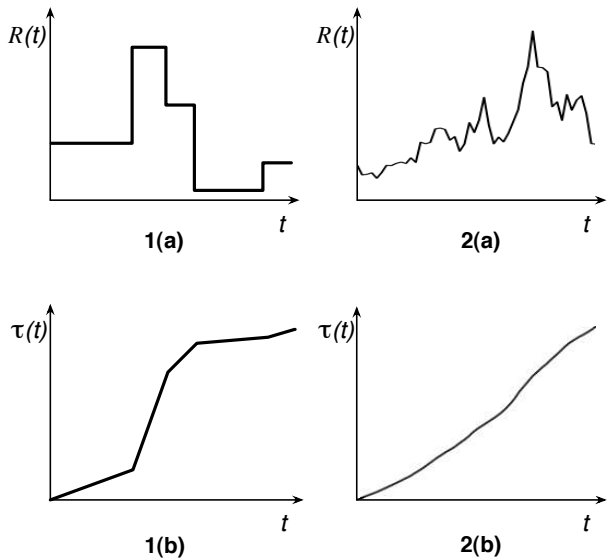}
\caption{\label{fig:discretecontinuous} A representation of the two classes of rate process with respect to the classification of property I. On top are examples of the rate history. Below are the corresponding integrated rates  $\tau(t) = \int_{s=0}^t R_s ds$. The figures 1(a-b) refer to a continuous-time Markov chain with discrete rate change events, and in figures 2(a-b), $R(t)$ is modelled as a diffusion process, with continuous paths.}
\vspace*{0.3in}
\end{center}
\end{figure}

The first characteristic is whether sample paths of the process are continuous or
discontinuous with respect to time.
Typically, models with discontinous paths have rates $R_t$ that
are constant except for discrete points in time at which there is
a jump in the value~(Figure~\ref{fig:discretecontinuous}-1(a)). If
the number of possible values for the rate is finite, then the
rate can easily be described as a continuous-time Markov chain
with a infinitesimal rate matrix. For example, in the covarion
process \cite{FM70} the basic rates are `off' ($R_t = 0$) or `on'
($R_t = 1$) and transitions occur between them at exponentially
distributed random time intervals. Galtier \cite{G01,G04}
generalizes this process to one with more than two possible
states. In other models, the range of possible values for the rate
is continuous, as in the model of Huelsenbeck \cite{HLS00}, where
a rate change event consists of multiplying the previous rate by a
gamma random variable. The rate change events are still discrete
and exponentially distributed.

There also are models that describes the rate as a continuous
function with time, and the most important class of Markov
processes with continuous paths are
\emph{diffusions}~(Figure~\ref{fig:discretecontinuous}-2(a)).
Examples include the CIR process presented here, the
Ornstein-Uhlenbeck model of Aris-Brosou and Yang \cite{AY02,AY03},
and the log-normal model of Kishino et al. \cite{KT01, TK98}.

Finally, it is also possible for $R_t$ to make jumps in value at a discrete set of times while also changing continuously in between these points.\\

\noindent{\em Property II: Long Term Behaviour and Ergodicity}\\

The second property we consider is the distribution of $R_t$ as
$t$ goes to infinity, that is, the distribution of the rate of
evolution in the long term. Surprisingly, many models of rate
evolution are very badly behaved in the limit.

One problematic class of processes that have already been applied
to rates in phylogenetics is the martingales. We say that a Markov
process is a \emph{martingale} if, for all $s,t \geq 0$ we have
$\EE [M_{t+s}| M_t] = M_t$~\cite{ok}. An example of of a Markov
martingale is Brownian motion. As a result of this fairly
innocuous looking condition, a martingale $M_t$ has the property
that either $\EE [|M_t|]$ is unbounded in time or $M_t$ converges
to a random constant \cite{Williams}. Either possibility is
undesirable from a modelling point of view. This may not be a
problem if we only look at the process over a finite time, but
neither is it particularly desirable. The processes of Kishino et
al.~\cite{KT01, TK98} and Huelsenbeck et al.~\cite{HLS00} all have
the property that either $R_t$ or $\log(R_t)$ is a martingale.

At a purely theoretical level, we observe that  an ever-increasing
variance will result for almost any signal that is \emph{only}
driven by its initial value and a stochastic force, with no
directional bias. The position of a particle subjected to a random
force produced by collisions with other particles is a classical
example of such a case. In our context, the effects on the
evolutionary rate are not independent of the actual rate :
whatever the theoretical framework we consider, a high
evolutionary rate is not as likely to increase (or to stay in high
values) as to go back to smaller values. The episodic evolution
fits particularly well to this idea, where periods of drastic
adaptation with high evolutionary rates are naturally followed by
periods where a population is adapted and its genome evolves much
more slowly. Even according to the neutral theory, as argued by
Takahata \cite{T87}, the overall dynamic of the rate should behave
as a random function that takes high values whenever bottlenecks
occur, and goes back to small values afterwards.

The concept of ergodicity naturally arises from this observation.
We say that a Markov process is \emph{ergodic} if for any initial
rate $R_0$ the distribution of $R_t$ converges to a unique
distribution as $t$ goes to infinity. The limiting distribution is
known as the invariant or stationary distribution. Examples of
ergodic processes include the OU process, the CIR model and
(usually) the discrete space covarion and covarion-type models
\cite{TS98, G01, G04}.

One possible way for a process to not be ergodic is if the
distribution of $R_t$ does not converge for large $t$ for some
initial rate $R_0$. This must be the case if $R_t$ is a martingale
and does not converge to a constant, as is the case with Brownian
motion. Another possibility is that $R_t$ converges to different
stationary distributions for different values of $R_0$.

\vspace{0.2in}

\noindent{\em Property III: Tractability} \\

A highly desirable feature of any model is its tractability, both
mathematical (does there exist a closed formula?) and
computational (can we compute probabilities efficiently?).
Nowadays, Monte Carlo methods make it possible to use arbitrarily
complex models: however, explicit analytical formulae allow for
more efficient sampling \cite{Liu}.

There are several probability distribution functions that are
important to have when working with rate processes. The most basic
is the distribution of the rate $R_t$ given the rate at time
$t=0$. This we have for the models \cite{AY02,AY03, TK98,KT01} and
for the CIR model, but not for the models of \cite{HLS00}.

In phylogenetics we incorporate the model for evolutionary rate into the mutation model for sequence evolution at a site. These interact to give a joint process $(R_t,X_t)$ for both the rate $R_t$ at time $t$ and the nucleotide or protein $X_t$ at time $t$. To evaluate the likelihood we require an expression for the joint conditional probability
\begin{eqnarray}
P[X_t = j, R_t =  s| X_0 = i, R_0 = r]
\end{eqnarray}
of going from one nucleotide (or amino acid) state and rate state
to another pair of states. Even though it is sometimes possible to
perform Monte Carlo computations to estimate this probability
without a formula (as in \cite{HLS00}), having a formula will
speed up the computations significantly without having to resort
to approximations,  as in \cite{TK98,KT01}.

\vspace{0.2in}

\noindent{\em Property IV: Autocovariance and dispersion}\\

There is general agreement \cite{G89}, \cite{CC93}, \cite{S97} on
the relevance of autocorrelation in the modelling of evolutionary
rate. Broadly speaking, if the various causes that explain rate
variation (generation time, population size, environmental
fitness) vary with time, it should be reflected in rate
variations. The extent to which the rate varies can be studied
using the {\em index of dispersion} (Kimura, \cite{K71}
\cite{K83}, Langley and Fitch \cite{F73}). Let $N(t)$ be the
number of substitutions or mutations of a sequence over time $t$.
The index of dispersion $I(t)$ is defined as
\begin{eqnarray}
I(t) = \frac{\mbox{Var}[N(t)]}{\EE [N(t)]}.
\end{eqnarray}
This statistic can be estimated by comparing the number of substitutions that have accumulated in different lineages  \cite{K83,B89}. The population genetics community has proposed different models to account for a high index of dispersion (\cite{T87}, \cite{C00B}), and any reasonable model should yield an index of dispersion of at least one.

The index of dispersion resulting from a particular model of rate variation is a function of the {\em autocovariance} of that model. The autocovariance for a process $R_t$ is defined by
\begin{eqnarray}
\rho(t) = \mbox{Cov}(R_0,R_t). 
\end{eqnarray}
For many processes we can derive an explicit formula for the autocovariance. If we assume that the substitutions occur according to a Poisson process with rate governed by our rate process (that is, the substitutions follow a {\em doubly stochastic} or {\em Cox} process, Section~\ref{MutModel}) and the rate process has autocovariance function $\rho(t)$ then
\begin{eqnarray}
I(t) = 1 + \frac{2 \int_0^t \left(1-\frac{s}{t} \right) \rho(s)  ds}{\EE (R_t)}, \label{IofD}
\end{eqnarray}
as stated by a theorem in \cite{CI80}, and the stationary index of dispersion \cite{G91} is then
\begin{eqnarray}
I(\infty) = \lim_{t \rightarrow \infty} I(t) = 1 + \frac{2 \int_0^{\infty} \rho(s)  ds}{\mu}, \label{IofDInf}
\end{eqnarray}
provided that $\mu$, the stationary mean of the process $R(t)$,
and the limit, exist. Note that if there is {\em any} variation in
rate then the index of dispersion will be greater than one
\cite{G91}.

Some rate models in phylogenetics \cite{TK98, AY03} don't model
explicitly the rate, but instead assign a (fixed) rate to each
branch, so that the expected number of substitutions on a
particular branch is equal to its length times its assigned
(constant) rate.

A close look at the log-normal model from Thorne et
al.~\cite{KT01}, which differs from their previous
version~\cite{TK98} in that the rate is explicitly modelled, we
suggest that the rate has constant autocovariance, since this
rate process is close to a transform of the Brownian motion, and Brownian motion 
has a constant autocovariance function. Put into
equation~(\ref{IofDInf}), we see that the the index of dispersion
diverges. This problematic result illustrates the necessity of a
balance between the presence of autocorrelation on one side, and
the decrease of autocorrelation on a large time scale.

\vspace{0.2in}

\noindent{\em Property V: Heterotachy or Homotachy}\\

There are two general ways that models for evolutionary rate can be incorporated into phylogenetics. On one hand, we have rate variation among lineages that applies to all sites (or loci)  together. This can be modelled by trees for which the paths from the root to the leaves have different lengths. The rate variation explains the extent to which the evolution of the sequences has violated the molecular clock.
Alternatively, we can introduce a distinct rate process for each site or locus.  This models {\em heterotachy}, where the lineage rate changes are site-specific \cite{LC02}. The transition probabilities that we derive in Section~\ref{MutModel} can be applied to homotachous as well as heterotachous models.


\begin{table}[h]
    \tiny
    \begin{center}
         \begin{tabular}{|p{1.6cm}p{1.5cm}|p{1.6cm}|p{1.4cm}|p{1.7cm}|p{1.7cm}|p{0.8cm}|}
            \hline
            & & I & II & III & IV & \\
            & & Rate Class &Ergodicity & Closed form for the transition probability & Closed form for the autocovariance & Ref.\\
            \hline
            \multirow{4}{1cm}{Models from population genetics} & Fluctuating mutation rates  & CTMC with continuous state space & No & Yes &None& \cite{T87} \\ 
            &  Fluctuating neutral space & CTMC with continuous state space &  Yes & None & None  & \cite{T87} \\ 
            &  Compound Poisson process   & CTMC with finite state space & Yes & None & None& \cite{T87, C00B} \\ 
            &  Episodic evolution  & CTMC with finite state space & Yes & None & None  &  \cite{G91}\\ \hline
            \multirow{5}{1cm}{Models from phylogenetics} & Covarion  & CTMC with finite state space & Yes &  Yes & None & \cite{FM70,TS98,G01,G04}  \\ 
            &  HLS   &  CTMC with continuous state space & No & None & None&  \cite{HLS00} \\ 
            &  Log-Normal   & Diffusion & No & None & Constant autocovariance & \cite{TK98,KT01}  \\ 
            &  Ornstein-Uhlenbeck & DIffusion & No & None  & Exponentially decreasing & \cite{AY02,AY03} \\ 
       &  CIR process & Diffusion & Yes & Yes  &  Exponentially decreasing&   \cite{CIR85} \\
        \hline
       \end{tabular}
            \caption{Models for the substitution rate, classified according to the properties of section~\ref{Prop} . CTMC stands for ``Continuous-Time Markov Chain".}  \label{ModelSummary}
       \end{center}
\end{table}


\section{A Continuous diffusion model for the evolutionary rate} \label{CIR}

A Markov process with continuous paths and satisfying some additional
smoothness conditions on its transition probabilities
\cite{KT81} is called a \emph{diffusion}. There are many
ways of specifying a diffusion process: perhaps the most intuitive one
is by giving the probability distribution function (pdf)  of $R_t$
given $R_0=r_0$, for arbitrary $r_0$.We denote this pdf by
$f_R[R_t|r_0]$. For example, Brownian motion with parameter $\sigma^2$
is defined by the condition that $f_R[R_t|r_0]$ is a normal
density with mean $r_0$ and variance $\sigma^2 t$.

A mathematically convenient representation of a diffusion is by
means of a \emph{stochastic differential equation} (SDE). In the same
way that a dynamical system can be defined as the solution of a
differential equation, a diffusion process $R_t$ can be defined as the
solution of an equation taking the general form (see \cite{ok} p.61)
\begin{eqnarray}
dR_t = \alpha(t,R_t)dt + \beta(t,R_t) dB_t. \label{eqn:sde}
\end{eqnarray}
Here, $\alpha(t,R_t)$ represents the deterministic effect on $R_t$,
$\beta(t,R_t)$ the stochastic part, and $dB_t$ is an infinitesimal
``random" increment. Brownian motion corresponds to the case when
$\alpha(t,R_t) = 0$ for all $t$, $\beta(t,R_t)$ is constant and the SDE becomes
\begin{eqnarray*}
dR_t =   \sigma dB(t).
\end{eqnarray*}
Note that if $\beta(t,R_t)=0$ for all $t$ and $R_t$ then \eqref{eqn:sde} becomes
a deterministic ordinary differential equation.

Going from an SDE such as \eqref{eqn:sde} to a pdf for the diffusion
involves solving a variable-coefficient second-order partial differential equation (PDE). For
general functions $\alpha$ and $\beta$ this PDE has no analytic solution.
There are very few diffusions known that
have closed form equations for their pdfs, and even fewer of these are
ergodic. The simplest ergodic diffusions with closed-form expressions
for the pdf are the
Ornstein-Uhlenbeck and the CIR (Cox-Ingersoll-Ross) \cite{CIR85}
processes.

The Ornstein Uhlenbeck (OU) process is described by the SDE
\[
dR_t = - b R_t dt +  \sigma dB_t.
   \]
The pdf for $R_t$ given $R_0 = r_0$ is the normal density
with mean $r_0 e^{-bt}$ and variance $\sigma^2 (1-e^{-2\theta
  t})$. Its stationary distribution is normal with mean 0 and
variance $\sigma^2$. The  OU process was used by Aris-Brosou and Yang
\cite{AY03} to model evolutionary rates. However, the OU process can take on negative values, and it is not clear how it can be used directly without any transformation, such as a reflected OU or a squared OU. Aris-Brosou and Yang also proposed another model, the EXP (for exponential) model, defined as the following : the rate assigned to a branch is drawn from an exponential distribution with mean equal to the rate of its ancestral branch. It is then obvious that their EXP model was a martingale. They outlined that the OU model seemed to provide a better fit to their data than the EXP model. Even if the reason of this better fit is still to be investigated, it seems reasonable to think that the ergodic property of the OU model could be a important factor. They also mentioned that the $\sigma^2$ parameter of the OU model was hard to infer, perhaps because the OU model has an insufficient number of free parameters.

The use of the CIR model solves the problem, since it is a generalization of the squared OU process, where the mean and the variance can be independently inferred by the addition of a third parameter. If the mean is fixed to one, we avoid any identifiability problem with branch lengths without fixing the variance, which can therefore be inferred as well as the autocorrelation.

The CIR process satisfies the SDE
\begin{equation}   \label{eqn:cirsde}
dR_t = b (a- R_t) dt + \sigma \sqrt{R_t} dB_t,
\end{equation}
and the pdf $f_R(R_t|r_0)$ for $R_t$ given $R_0 = r_0$ is a non-central $\chi^2$
distribution with degree of freedom $4ab/\sigma^2$ and parameter of non-centrality $\frac{4 b r_0e^{-bt}}{\sigma^2 (1-e^{-bt})}$. Its mean and variance are equal to

\begin{eqnarray}
\EE[r_t] &=& r_0 e^{-b t} + a(1-e^{-bt}) \\
\mathrm{Var}[r_t] &=&  r_0 \frac{\sigma^2}{b} (e^{-bt} - e^{-2bt}) + \frac{a \sigma^2}{2b}(1-e^{-bt})^2.
\end{eqnarray}

The stationary distribution of $R_t$ is
a gamma distribution with shape parameter $2ab/ \sigma^2$ and  scale
parameter $\sigma^2/2b$. Hence the mean of the stationary distribution
is  $a$ and the variance is $\frac{a \sigma^2}{2b}$ \cite{CIR85}.

Unlike an OU process, if $r_0$, $a$, and $b$ are all positive a CIR process is always non-negative. The square
of an OU process is a special case of the CIR process. Furthermore, by multiplying $R_t$ by a constant in equation~(\ref{eqn:cirsde}), we see that multiplying a CIR process by a positive constant gives another CIR process.

The covariance of the stationary CIR process can be exactly computed as
\begin{eqnarray}
\rho(t) = \mbox{Cov}(R_0,R_t) = \frac{a\sigma^2 }{2b} e^{-bt}.
\end{eqnarray}
From this, (\ref{IofD}) leads to a closed formula for the index of dispersion:
\begin{eqnarray*}
I_{CIR}(t) = 1 + \frac{\sigma^2}{b^3 t}  (bt -1 + e^{-bt}).
\end{eqnarray*}
Thus
\begin{equation}
I_{CIR}(\infty) = \lim_{t \rightarrow \infty} I_{CIR}(t) =  1+ \frac{\sigma^2}{b^2} . \label{eqn:Icir}
\end{equation}

In Section~\ref{Prop} we emphasized that the concept of autocovariance is close to the index of dispersion. As Zheng showed~\cite{Z01}, the effect of complex infinitesimal rate matrices on the index of dispersion (with constant rate) is not likely to explain alone the observed large empirical values. If the rate varies, Cutler~\cite{C00C} showed that an elevated index of dispersion can only be achieved if the time-scale of the rate process is approximately of the same magnitude as the substitution process itself. The CIR process provides the possibility to satisfy this property, while incorporating autocovariance. It is the consensus of these two ideas that should guide our choice for the rate of evolution.

From \eqref{eqn:cirsde} we see that the CIR process possesses three parameters, namely the stationary mean $a$, the stationary variance $\sigma^2$, and the intensity of the force that drives the process to its stationary distribution, $b$. The parameter $b$ determines how fast the process autocovariance goes to 0 as $t$ increases.

The three parameters of the CIR process can be quickly estimated from standard statistics in molecular evolution. The parameter $a$ is a scale parameter. It determines the expected rate at any time  given no other information. Throughout the paper, we will assume that $a=1$, so that the model has an expected rate equal to one. This parallels the constraint that the gamma distribution has an expected rate equal to one in the Rate-Across-Site (RAS) model \cite{Y93}.

The CIR process has a stationary distribution given by a gamma distribution. To make the stationary distribution coincide with the gamma distribution of a RAS model with parameter $\Gamma$ we choose $\sigma$ and $b$ such that
\begin{equation} \Gamma = \frac{\sigma^2}{b}. \label{eqn:gammaform} \end{equation}
The stationary index of dispersion, $I_{CIR}(\infty)$, can be estimated empirically \cite{F73, G89}. We can then use \eqref{eqn:Icir} and \eqref{eqn:gammaform} to obtain the estimates
\begin{eqnarray*}
\hat{b} &=& \frac{\hat{\Gamma}}{\hat{I}_{CIR}(\infty) -1} ,\\
\hat{\sigma}^2 &=& \frac{\hat{\Gamma}^2}{\hat{I}_{CIR}(\infty)-1}.
\end{eqnarray*}

\section{Mutation models with a rate process} \label{MutModel}

The standard model for the substitution process  at a particular locus is a continuous-time Markov chain. This kind of process is defined by a square matrix $Q$ called the \emph{infinitesimal rate matrix}. Suppose, to begin, that there is a constant evolutionary rate $r_0$. As above, we let $X_t$ denote the state (e.g. amino acid) at time $t$. The transition probabilities are then given by
\begin{equation}
\Pr[X_t = j|X_0 = i] = [e^{Qr_0t}]_{ij}. \label{eqn:transprobsConstant}
\end{equation}
We suppose that the process has a unique stationary distribution $\pi$, where $\pi_j$ is the stationary probability of state $j$ and
\[\pi_j = \lim_{t \rightarrow \infty} \Pr[X_t = j|X_0 = i]\]
for all $i$ and $j$.  We assume that $Q$ has been normalised so that in the stationary distribution the expected number of mutations over time $t$ equals $r_0t$.
Note that the transition probabilities \eqref{eqn:transprobsConstant} depend only on the product $r_0t$, so will be the same if we double the rate and halve the time, for example.

Suppose now that the rate is not constant, but instead varies according to some fixed function $r_s$, $s \geq 0$. Equation \eqref{eqn:transprobsConstant} then becomes
\begin{equation}
\Pr[X_t = j|X_0 = i, r] = [e^{Q \tau_r}]_{ij}. \label{eqn:transprobsFunction}
\end{equation}
where
\[\tau_r = \int_{s = 0}^{s=t} r_s ds\]
is the area under the curve $r_s$.

In the models we will consider, the fixed function $r = (r_t)_{t \geq 0}$ is replaced by a random process $R = (R_t)_{t \geq 0}$ that is dependent only on the starting rate $r_0$. The integral
\begin{equation}
\tau_R = \int_{s=0}^{s=t} R_s ds \label{eqn:taur}
\end{equation}
is also random in this case; let $g_{\sR}$ denote its pdf. The transition probabilities  can be determined from the expected value of \eqref{eqn:transprobsFunction} with $\tau_r$ replaced by the random variable $\tau_{\sR}$. By the law of total expectation, this simplifies to
\begin{equation}
\Pr[X_t = j|X_0 = i] = \int_{\tau} [e^{Q \tau}]_{ij} g_{\sR} (\tau) d\tau  \label{eqn:transprobsProcess}.
\end{equation}

Let $M(\eta) = \EE_\tau[e^{\eta \tau_{\sR}}]$ denote the {\em moment generating function} (mgf) for the random variable $\tau_{sR}$. Then \eqref{eqn:transprobsProcess} can be rewritten
\[\Pr[X_t = j|X_0 = i] = [M(Q)]_{ij}\]
where the function $M$ is interpreted as a matrix function \cite{HJ85}. We assume that $Q$ can be diagonalised as $Q =  V \Lambda V^{-1}$, where  $\Lambda =  \mbox{diag}(\lambda_1, \ldots, \lambda_n)$ is a diagonal matrix formed from the eigenvalues of $Q$. The matrix function $M(Q)$ can then be evaluated as $M(Q) = V M(\Lambda) V^{-1}$, where
\[M(\Lambda) =  \mbox{diag}(M(\lambda_1), \ldots, M(\lambda_n)).\]
See \cite{matrix}  for a more details on matrix functions. The problem of determining pattern probabilities therefore boils down to the problem of determining the moment generating function of the integrated rate, $\tau_{\sR}$ (eqn. \ref{eqn:taur}). Tuffley and Steel use this approach to derive distance estimates for the covarion process \cite{TS98}. 

For applications in phylogenetics, we need the mgf of $\tau_{\sR}$ conditioned on just the starting rate, or both the starting and finishing rate. The mgf of $\tau_R$, conditioned on a starting rate of $r_0$, is then
\begin{eqnarray}
M_{r_0}(\eta) & = & \EE\big[\exp(\eta \tau_R) | R_0 = r_0 \big] \nonumber \\
& = & \EE\left[\exp\left(\eta \int_{s=0}^t R_s ds \right) \Big| R_0 = r_0 \right]. \label{eqn:Mr0}
\end{eqnarray}
As before, we let $f_R(R_t|r_0)$ denote the pdf of $R_t$ conditioned on $R_0 = r_0$. Let $\delta(x)$ denote the Dirac delta function with $\delta(0) = 1$ and $\delta(x) = 0$ for all $x \neq 0$. The 
mgf of $\tau_R$ conditioned on both the starting and finishing rates is
\begin{eqnarray}
\hspace*{-0.5cm} M_{r_0,r_t}(\eta) & = & \EE\big[\exp(\eta \tau_R) \big| R_0 = r_0,\,R_t = r_t\big] \nonumber\\
& = & \frac{1}{f_R(r_t|r_0)} \EE\big[\exp(\eta \tau_R) \delta(R_t - r_t) \big| R_0 = r_0\big] \nonumber\\ 
& = &  \frac{1}{f_R(r_t|r_0)}  \EE\left[\exp\left(\eta \int_{s=0}^t R_s ds \right)  \delta(R_t - r_t) \,\,\Big | \,\, R_0 = r_0\right]. \label{eqn:Mr0rt}
\end{eqnarray}
Equations \eqref{eqn:Mr0} and \eqref{eqn:Mr0rt} hold irrespective of whether $R$ is discrete or continuous, a diffusion, jump process, or a continuous time Markov chain. 

We note in passing that analytic formulae for $M_{r_0}(\eta)$ and $M_{r_0,r_t}(\eta)$ exist in the case that $R$ is a continuous time Markov chain, for example in the covarion-type model of Galtier  \cite{G01}.  Suppose that the evolutionary rate switches between rate values $g_1,g_2,\ldots g_k$ following a continuous time Markov chain with infinitesimal rate matrix $G$. Let $D$ be the $k \times k$ diagonal matrix with entries $g_1,g_2,\ldots,g_k$. A careful reworking of the proof of Theorem 1 in \cite{DM68} gives the mgf of $\tau_{\sR}$ conditioned on both the starting and finishing rate.  The mgf for $\tau_R$ conditioned on $r_0 = g_i$ is then
\[M_{g_i} =  \sum_{j=1}^k (e^{(G+\eta D)t})_{ij} \]
while the mgf of $\tau_R$ conditioned on $r_0 = g_i$ and $r_t = g_j$ is
\[M_{g_i,g_j} = \frac{(e^{(G+\eta D)t})_{ij}}{(e^{Gt})_{ij}} .\]
This provides an independent derivation of the formula in \cite{G04} for transition probabilities under a covarion-type model.

\section{Moment generating functions and transition probabilities for the CIR model} \label{Sectionmgf}

In this section we derive expressions for the (joint) transition probabilities
\begin{eqnarray}
\Pr[X_t = j | X_0 = i, R_0 = r_0]. \label{jointpdfStart}
\end{eqnarray}
and
\begin{eqnarray}
\Pr[X_t = j, R_t =  s| X_0 = i, R_0 = r_0]\label{jointpdf}
\end{eqnarray}
As we have seen, to evaluate these probabilities we need to determine the moment generating functions (mgfs)  defined in equations \eqref{eqn:Mr0} and \eqref{eqn:Mr0rt}.

We use the Feynman-Kac formula \cite{Kac,ok} to derive analytic formulae for $M_{r_0}(\eta)$ and $M_{r_t,r_0}(\eta)$ under the CIR model. Let $g(\cdot)$ be a real-valued function. Define the function $v = v(t,x)$ by
\begin{eqnarray}
v(t,x) = \EE \left[ \exp\left(  \eta \int_0^t R(s) ds \right) g(R_t)  \,\,\Big | \,\, R_0 = x\right]  \label{eqn:FK}
\end{eqnarray}
The Feynman-Kac formula \cite{Kac} asserts
 that $v(t,x)$ solves the following partial differential equation (PDE)
 \begin{eqnarray} \label{PDE}
 \frac{\partial }{\partial t} v(t,x)=  b(1-x) \frac{\partial}{\partial x} v(t,x) + \frac{1}{2} \sigma^2 x(t,x) \frac{\partial^2 }{\partial x^2} v(t,x) + \eta xv
  \end{eqnarray}
for $t>0$, $x \in \mathbb{R}$, and with boundary condition
\begin{eqnarray}
v(0,x) = g(x)\,\,\,\, \mbox{for all $x \in \mathbb{R}$.} \label{eqn:boundary}
\end{eqnarray}
We apply the methods in   \cite{S04} and \cite{AL04} to solve these pdes with the different boundary conditions.

First consider the case when we condition only on the initial rate, eqn.~\eqref{eqn:Mr0}. To make \eqref{eqn:FK} equal to \eqref{eqn:Mr0} we set $g(x) = 1$ for all $x$.  The boundary condition \eqref{eqn:boundary} in this case becomes 
\[ v(0,x) = 1\,\,\,\, \mbox{for all $x \in \mathbb{R}$.} .\]
With this boundary condition, the pde \eqref{PDE} has solution
\[v(t,x) =  \Psi(\eta,t) e^{-x \Xi(\eta,t)}\]
where
\begin{eqnarray}
\Psi(\eta,t) & =& \left(  \frac{\overline{b} e^{bt/2}}{\overline{b} \cosh(\overline{b} t/2) + b \sinh(\overline{b}t/2)} \right)^{\frac{2b}{\sigma^2}}\\
\Xi(\eta,t) &=& \left( \frac{2 \eta \sinh(\overline{b}t/2)}{\overline{b}\cosh(\overline{b}t/2) + b \sinh(\overline{b}t/2)} \right), \label{Xi} \\
\overline{b} &=& \sqrt{b^2 - 2 \eta \sigma^2}
\end{eqnarray}
We therefore have
\begin{eqnarray}
M_{r_0}(\eta) = \Psi(\eta,t) e^{-r_0 \Xi(\eta,t)}. \label{mgfCIRcond}
\end{eqnarray}

The case when both the starting and finished rates are specified is more complicated. From \eqref{eqn:Mr0rt} the mgf $M_{r_0,r_t}$ can be written
\[ M_{r_0,r_t} = \frac{1}{f_R(r_t|r_0)}  v(t,x)\]
where, in this case, $v(t,x)$ is given by \eqref{eqn:FK} with $g(x) = \delta(x - r_t)$. The boundary condition  \eqref{eqn:boundary} therefore becomes 
\[ v(0,x) = \delta(x - r_t) .\]
With this new boundary condition, the pde \eqref{PDE} has solution
\begin{eqnarray}
v(t,x) &=&  c \exp \left[ -\frac{bt}{\sigma^2}(\overline{b}-b) + \frac{b-\overline{b}}{\sigma^2} x - \frac{b +\overline{b}}{\sigma^2} r_t- c(r_t + x) e^{-\overline{b}t} \right] \nonumber \\
& & \quad \quad \quad \times \left( \frac{r_t}{x e^{-\overline{b}t}} \right)^{\frac{b}{\sigma^2}-1/2} I_{\frac{2b}{\sigma^2}-1} \left( 2c\sqrt{x r_t e^{-\overline{b}t}} \right), \label{solveCIR}
\end{eqnarray}
where
\begin{eqnarray*}
c & =& \frac{2\overline{b}}{\sigma^2(1-e^{-\overline{b}t})} \\
\overline{b} &=& \sqrt{b^2 - 2 \eta \sigma^2}
\end{eqnarray*}
and $I_{\nu}(x)$ is the modified Bessel function of the first kind with parameter $\nu$ \cite{AS}.  
Hence the mgf conditioned on initial and final rates is given by
\begin{eqnarray*}
M_{r_0,r_t}(\eta) &=&  c \exp \left[ -\frac{bt}{\sigma^2}(\overline{b}-b) + \frac{b-\overline{b}}{\sigma^2} r_0 - \frac{b +\overline{b}}{\sigma^2} r_t- c(r_t + r_0) e^{-\overline{b}t} \right]  \\
& &   \quad \times \left( \frac{r_t}{r_0 e^{-\overline{b}t}} \right)^{\frac{b}{\sigma^2}-1/2} I_{\frac{2b}{\sigma^2}-1} \left( 2c\sqrt{r_0 r_t e^{-\overline{b}t}} \right) \frac{1}{f_R(r_t|R_0=r_0)}, \end{eqnarray*}
where $c$ and $\overline{b}$ are defined above and, from Section~\ref{CIR}, $f_R(R_t|R_0 = r_0)$  is the pdf for a non-central $\chi^2$
distribution with degree of freedom $4ab/\sigma^2$ and parameter of non-centrality $\frac{4 b r_0e^{-bt}}{\sigma^2 (1-e^{-bt})}$,

Bringing everything together, we have our main result.

\begin{theorem}
Define $P$ by $P_{ij} = \Pr[X_t = j, R_t = r_t | X_0 = i, R_0 = r_0]$. Suppose that $Q = V \Lambda V^{-1}$ where $\Lambda$ is a diagonal matrix containing the eigenvalues $\lambda_1,\ldots,\lambda_n$ of $Q$. Then
\begin{eqnarray}
P = M(Q) = V M(\Lambda) V^{-1} \label{Transmgf}
\end{eqnarray}
where $M(\Lambda)$ is the diagonal matrix where, for all $i$,
\[M(\Lambda)_{ii} = M_{r_t,r_0}(\lambda_i).  \]
\end{theorem}

\section{Three taxa phylogenies} \label{sec:3taxa}

The simplest phylogeny for which we can distinguish between constant and variable evolutionary rates is a tree with three taxa. For this reason, there are many methods for testing, and estimating, rate variation that are based on three taxon analyses \cite{G89}. Here we show that the likelihood for a three taxa tree, under the CIR model of rate variation, can be computed exactly. The problem for general phylogenies is more complex since we have to integrate out rates for the internal nodes. Here, we consider a heterotachous model, so that each site has its own rate history. Because the sites (and the rate at each site) evolve independently from each other, the likelihood of a sequence will be the product of all site-specific likelihoods. Therefore, we only require the likelihood computation for one site.

We recall that the stationary distribution of the CIR is a gamma distribution $\Gamma(\omega,\nu)$, where $\omega = 2b/\sigma^2$ and $\nu=2b/\sigma^2$, i.e.
\begin{eqnarray}
f_{R_0}(r) = \frac{\omega^{\nu}}{\Gamma(\nu)} r^{\nu-1}e^{-\omega r} \label{statDist}.
\end{eqnarray}
Therefore the stationary mean and variance are $1$ and $\sigma^2/2b$.

In order to get the transition probabilities, we will use the mgf of $\tau_{\sR}$ unconditioned on the final rate, given by equation (\ref{mgfCIRcond}). The transition probability matrix of the subsitution process, given initial rates, can be obtained by equations (\ref{mgfCIRcond}) and (\ref{Transmgf}). Let $\lambda_1,\ldots,\lambda_n$ be the eigenvalues of $Q$. Using eigenvalue decomposition,
we can find vectors $\bu^{(1)},\ldots,\bu^{(n)}$ and $\bv^{(1)},\ldots,\bv^{(n)}$ such that
\begin{equation}
\Pr[X_t = i | X_0 =j, R_0 = r_0] = \sum_{k=1}^n {\bu^{(k)}_{j}}^{T} \bv^{(k)}_{i} M_{r_0}(\lambda_k,t), \label{eqn:trans}
\end{equation}
where we changed slightly our notation and explicitly wrote the dependency of $M_{r_0}$ on $t$.

\begin{figure}[h]
\vspace*{0.3in}
\begin{center}
\includegraphics[scale = 0.5]{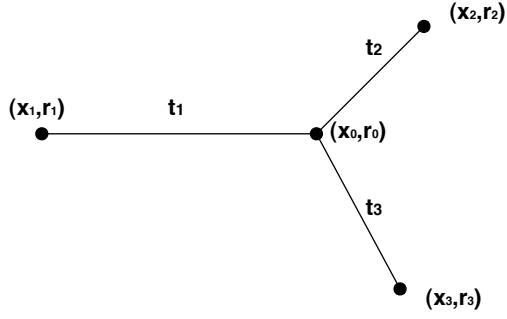}
\caption{\label{3taxatree} A three taxa unrooted tree, with branch lengths and one character state and rate value associated to each leaf.}
\vspace*{0.3in}
\end{center}
\end{figure}

Now consider the 3-taxa tree with branches of lengths $t_1,t_2,t_3$ leading to leaves labelled with states $x_1,x_2, x_3$ (Figure~\ref{3taxatree}). If we condition on a rate $r_0$ and state $x_0$ at the root
then the probability of observing $x_1,x_2,x_3$ at the leaves is given by
\begin{eqnarray}
L(x_1,x_2,x_3|x_0, r_0)&=& P[X_{t_1} = x_{1} | x_0,  r_0] P[X_{t_2} = x_2 | x_0,  r_0]  P[X_{t_3} = x_3 | x_0,  r_0] \nonumber \\
&=& \sum_{i=1}^n \sum_{j=1}^n \sum_{k=1}^n B_{ijk} M_{r_0}(\lambda_i,t_1) M_{r_0}(\lambda_j,t_2) M_{r_0}(\lambda_k,t_3) \label{eqn:Pexpansion}
\end{eqnarray}
where
\[B_{ijk} =  \bu^{(i)}_{x_0} \bv^{(i)}_{x_1} \bu^{(j)}_{x_0} \bv^{(j)}_{x_2} \bu^{(k)}_{x_0} \bv^{(k)}_{x_3}.\]

The rate at the root is assumed to have the stationary distribution $f_{R_0}$ given by \eqref{statDist}.
The likelihood integrated with respect to $r_0$ is then
\begin{eqnarray}
L(x_1,x_2,x_3|x_0) &=&  \int_{r_0} L(x_1,x_2,x_3|x_0, r_0) f_{R_0}(r_0) dr_0 \nonumber
\end{eqnarray}
which by \eqref{eqn:Pexpansion} equals
\begin{equation}
  \sum_{i=1}^n \sum_{j=1}^n \sum_{k=1}^n B_{ijk}  \int_{r_0} M_{r_0}(\lambda_i,t_1) M_{r_0}(\lambda_j,t_2) M_{r_0}(\lambda_k,t_3)  f_{R_0}(r_0) dr_0. \label{eqn:intlike}
\end{equation}
We now use the formula \eqref{mgfCIRcond} for the mgf's derived above.
\begin{eqnarray}
\lefteqn{M_{r_0}(\lambda_i,t_1) M_{r_0}(\lambda_j,t_2) M_{r_0}(\lambda_k,t_3)  f_{R_0}(r_0) dr_0} \nonumber \\
&=& \Psi(\lambda_{i},t_1) e^{-r_0 \Xi(\lambda_{i},t_1)} \Psi(\lambda_{j},t_2) e^{-r_0 \Xi(\lambda_{j},t_2)} \Psi(\lambda_{k},t_3) e^{-r_0 \Xi(\lambda_{k},t_3)}  \frac{\omega^{\nu}}{\Gamma(\nu)} r_0^{\nu-1}e^{-\omega r_0} \nonumber \\
&=& \Psi(\lambda_{i},t_1) \Psi(\lambda_{j},t_2) \Psi(\lambda_{k},t_3)  \frac{\omega^{\nu}}{\Gamma(\nu)} r_0^{\nu-1} e^{-r_0 \left( \omega + \Xi(\lambda_{i},t_1) +  \Xi(\lambda_{j},t_2) +  \Xi(\lambda_{k},t_3) \right)} \label{eqn:4'}
\end{eqnarray}
Using integration by parts, or simply using the fact that the gamma distribution integrates to 1, we get%
\begin{eqnarray}
\int_{r_0} \lefteqn{M_{r_0}(\lambda_i,t_1) M_{r_0}(\lambda_j,t_2) M_{r_0}(\lambda_k,t_3)  f_{R_0}(r_0) dr_0} \nonumber \\
&&\hspace*{-0.6cm}= \,\,\Psi(\lambda_{i},t_1) \Psi(\lambda_{j},t_2) \Psi(\lambda_{k},t_3) \left(\frac{\omega}{\omega + \Xi(\lambda_{i},t_1) +  \Xi(\lambda_{j},t_3) +  \Xi(\lambda_{k},t_3)} \right)^{\nu}. \nonumber
\end{eqnarray}
Finally, we can substitute this back into \eqref{eqn:intlike} to obtain
\begin{eqnarray}
L(x_1,x_2,x_3|x_0) &=&  \sum_{i=1}^n \sum_{j=1}^n \sum_{k=1}^n B_{ijk}  \Psi(\lambda_{i},t_1) \Psi(\lambda_{j},t_2) \Psi(\lambda_{k},t_3) \nonumber \\
&&\hspace*{1cm}\times \left(\frac{\omega}{\omega + \Xi(\lambda_{i},t_1) +  \Xi(\lambda_{j},t_3) + \Xi(\lambda_{k},t_3)} \right)^{\nu}. \nonumber
\end{eqnarray}

The formula extends immediately to phylogenies with $n$ leaves attached to the root, though the number of terms in the summation increases exponentially. Our approach has been to evaluate likelihoods on complete phylogenies using Monte-Carlo techniques, together with the exact transition probabilities derived here.

\section{Discussion}

\subsection{Summary}

We have shown how, given a few natural criteria for our model selection, the CIR appeared as the simplest continuous model that is at the same time ergodic, has a non-zero autocovariance function and that can account for an arbitrarily large index of dispersion. Moreover, we provided simple ways to estimate its parameters with the help of two observable statistics, namely the RAS gamma parameter and empirical index of dispersion. Another very interesting practical aspect of the CIR process is that it can be easily, and without approximations, implemented in the MCMC framework.

\subsection{Future extensions}

A possible future extension of our model could involve jump models, in which the rate path is discontinuous as in the continuous-time Markov chain, but also varies as diffusion between these discontinuities. However, the use of such a model implies the use of more parameters, and it may well be the case that the relative weakness of the rate of evolution signal cannot allow the use of more than two parameters, because of identifiability problems.


\begin{thebibliography}{999}

\bibitem{Y93}
Z. Yang, Maximum likelihood estimation of phylogeny from DNA sequences when substitution rates differ over sites,  Mol. Biol. Evol. 10 (1993) 1396-1401.

\bibitem{LC02}
P. Lopez, D. Casane, H.  Philippe, Heterotachy, an important process of protein evolution, Mol. Biol. Evol. 19 (2002) 1-7.

\bibitem{SI02}
E. Susko, Y. Inagaki, C. Field, M. E. Holder, A. J. Roger,  Testing for differences in rates-across-sites distributions in phylogenetic subtrees, Mol. Biol. Evol. 19 (2002) 1514-1523.

\bibitem{AY02}
S. Aris-Brosou, Z. Yang, Effects of models of rate evolution on estimation of divergence dates with special reference to the metazoan 18S ribosomal RNA phylogeny, Syst. Biol. 51 (2002) 703-714.

\bibitem{KT01}
H. Kishino, J. L. Thorne, W. J. Bruno, Performance of a divergence time estimation method under a probabilistic model of rate evolution, Mol. Biol. Evol. 18 (2001) 352-361.

\bibitem{B86}
R. J. Britten, Rates of DNA sequence evolution differ between taxonomic groups, Science 231 (1986) 1393-1398.

\bibitem{W89}
K. H. Wolfe, P. M. Sharp, W.-H. Li, Mutation rates differ among regions of the mammalian genome, Nature 337 (1989) 283-285.

\bibitem{MP93}
A. P. Martin, S. R. Palumbi, Body size, metabolic rate, generation time, and the molecular clock, Proc. Natl. Sci. USA 90 (1993) 4087-4091.

\bibitem{GY94}
N. Goldman, Z. Yang, A codon-based model of nucleotide substitution for protein-coding DNA sequences, Mol. Biol. Evol. 11 (1994) 725-736.

\bibitem{MG94}
S.V. Muse, B.S.  Gaut, A likelihood method for comparing synonymous and nonsynonymous nucleotide substitution rates, with application to the chloroplast genome, Mol. Biol. Evol. 11 (1994) 715-724.

\bibitem{PPH01}
C. P\`{a}l, B. Papp, L. D. Hurst, Highly expressed genes in yeast evolve slowly, Genetics 158 (1998): 927-931.

\bibitem{RJKGT03}
D. M. Robinson, D. T. Jones, H. Kishino, N. Goldman, J.L. Thorne, Protein evolution with dependence among codons due to tertiary structure, Mol. Biol. Evol. 20 (1998) 1692-1704.

\bibitem{OT90}
T. Ohta, H.  Tachida, Theoretical study of near neutrality. I. Heterozygosity and rate of mutant substitution, Genetics 126 (1990) 210-229.

\bibitem{T87}
N. Takahata, On the overdispersed molecular clock, Genetics 116 (1987) 169-179.

\bibitem{L69}
C. D. Laird, B. L. Mc Conaughy, B. J. Mc Carthy, Rate of fixation of nucleotide substitutions in evolution,  Nature 224 (1969) 149-154.

\bibitem{CIR85}
J. C. Cox, J. E.  Ingersoll, S. A. Ross, A theory of the term structure of interest rates, Econometrica 53 (1985) 385-408

\bibitem{UC71}
T. Uzzell, K.W. Corbin, Fitting discrete probability distributions to evolutionary events, Science 172 (1971) 1089-1096.

\bibitem{FM70}
W. M. Fitch, E. Markowitz, An improved method for determining codon variability in a gene and its application to the rate of fixation of mutations in evolution, Biochem. Genet. 4 (1970) 579-593.

\bibitem{G01}
N. Galtier, Maximum-likelihood phylogenetic analysis under a covarion-like model, Mol. Biol. Evol. 18 (2001) 866-873.

\bibitem{G04}
N. Galtier, Markov-modulated markov chains and the covarion process of molecular evolution, Journal of Computational Biology 11 (2004) 727-733.

\bibitem{HLS00}
J. P. Huelsenbeck, B. Larget, D. Swofford, A compound Poisson process for relaxing the molecular clock, Genetics 154 (2000) 1879-1892

\bibitem{AY03}
S. Aris-Brosou, Z. Yang, Bayesian models of episodic evolution support a late precambrian explosive diversification of the metazoa, Mol. Biol. Evol. 20 (2003) 1947-1954.

\bibitem{TK98}
J.L. Thorne, H. Kishino, I.S. Painter, Estimating the rate of evolution of the rate of molecular evolution, Mol. Biol. Evol. 15 (1998) 1647-1657.

\bibitem{ok}
B. \O ksendal, Stochastic Differential Equations: an Introduction with Applications, Springer, 1998.

\bibitem{Williams}
D. Williams, Probabilities with Martingales, Cambridge mathematical textbooks, 1991.

\bibitem{TS98}
C. Tuffley, M. A. Steel, Modeling the covarion hypothesis of nucleotide substitution, Math. Biosci. 147 (1998) 63-91.

\bibitem{Liu}
J. S. Liu, Monte Carlo Strategies in Scientific Computing, Springer, 2001.

\bibitem{G89}
J.H. Gillespie, Lineage effects and the index of dispersion of molecular evolution, Mol. Biol. Evol. 6 (1989) 636-647.

\bibitem{CC93}
L. Chao, D. E. Carr, The molecular clock and the relationship between population size and generation time, Evolution 47 (1993) 688-690.

\bibitem{S97}
M. J. Sanderson, A  nonparametric approach to estimating divergence times in the absence of rate constancy, Mol. Biol. Evol. 14 (1997) 1218-1231.

\bibitem{K71}
T. Ohta, M. Kimura, On the constancy of the evolutionary rate in cistrons, J. Mol. Evol. 1 (1971) 18-25.

\bibitem{K83}
M. Kimura, The Neutral Theory of Molecular Evolution, Cambridge University Press, Cambridge, 1983.

\bibitem{F73}
C.H. Langley, W.M. Fitch, The constancy of evolution: a statistical anlalysis of the $\alpha$ and $\beta$ haemoglobins, cytochrome c, and fibrinopeptide A, in: Genetic Structure of Populations, Univ. of Hawaii press, Honolulu, 1973

\bibitem{B89}
M. Bulmer, Estimating the variablility of substitution rates, Genetics 123 (1989) 615-619.

\bibitem{C00B}
D. J. Cutler, Understanding the overdispersed molecular clock, Genetics 154 (2000) 1403-1417.

\bibitem{CI80}
D. R. Cox, V. Isham, Point Processes, Chapman and Hall, NewYork, 1980.

\bibitem{G91}
J. H., Gillespie, The Causes of Molecular Evolution, Oxford University Press, 1991.

\bibitem{KT81}
S. Karlin, H. M. Taylor, H. M., A Second Course in Stochastic Processes, New York academic press, 1981.

\bibitem{Z01}
Q. Zheng, On the dispersion index of a Markovian molecular clock, Math. Biosci. 172 (2001) 115-128.

\bibitem{C00C}
D. J. Cutler, The index of dispersion of molecular evolution: slow fluctuations, Theor. Pop. Biol. 57 (2000) 177-186.

\bibitem{HJ85}
R. Horn, C. Johnson, Matrix Analysis, Cambridge University Press, 1985.

\bibitem{matrix}
G. H. Golub, C. F. Van Loan, Matrix Computations, John Hopkins University Press, 1996.

\bibitem{DM68}
J. N. Darroch, K. W. Morris, Passage-time generating functions for continuous-time finite Markov chains, J. Appl. Prob. 5 (1968) 414-426.

\bibitem{Kac}
M. Kac, On Some connections between probability theory and differential and integral equations, Proceedings of the second Berkeley symposium on probability and statistics, University of California, Berkeley, 1951.

\bibitem{S04}
S. E. Shreve, Stochastic Calculus for Finance, Springer Finance, 2004.

\bibitem{AL04}
C. Albanese, S. Lawi, Laplace Transforms for Integrals of Markov Processes, Submitted to Markov Proc. Rel. Fields (2005).

\bibitem{AS}
M. Abramowitz, I. A. Stegun, Handbook of Mathematical Functions with Formulas, Graphs, and Mathematical Tables, Dover New York, 1965.







\end{thebibliography}
\end{document}